\documentclass{ifacconf}
\usepackage{graphicx}      
\usepackage{natbib}        
\usepackage{newfloat}
\usepackage{listings}
\usepackage{amsmath}
\usepackage{amsfonts}
\usepackage{algorithm}
\usepackage{algorithmic}

\begin{document}
\begin{frontmatter}
\title{Dynamics and Stability under Iterated Sanctions and Counter-sanctions}
\author{Alexander Fradkov}
\address{Institute for Problems in Mechanical Engineering, \\ 61, Bolshoy pr .V.O., Saint Petersburg 199178 Russia \\}
\address{St Petersburg University \\ 28 Universitetskii prospect, Peterhof,
St Petersburg 195904, Russia}
\date{July 2022}

\begin{abstract}
In the paper a simple discrete dynamical model for dynamics of two antagonistic agents (opponents) under iterated sanctions and counter-sanctions is proposed.  The model is inspired with Osipov-Lanchester model for combats. Simple stability criteria are derived both for the full information case and for stochastic uncertainty case. The results provide some important qualitative conclusions that are interpreted in terms of international stability preservation.

\end{abstract}

\end{frontmatter}
\section{Introduction}
\par Imposing sanctions is a traditional tool of both international and internal affairs. The classic models of military operations and battles - the Osipov-Lanchester model \cite{Osipov1915,Lanchester1916, Helmbold93} and its modifications are still being studied in various educational institutions \cite{McCartney22,Johnson15,Kress20,Washburn09}. However, the processes of antagonistic interaction between  competing agents or opponents are not limited to military confrontation. They also include economic and political actions called sanctions. The term "sanctions" has become often used and even popular during recent years, but in fact, the processes of economic and political confrontation in the modern world existed many decades ago.
In 2021-2022, sanctions began to be used especially often, and sometimes even began to determine the evolution of international relations.

The role of sanctions is similar to the role of military operations. However, mathematical models of sanctions differ from the models of military operations and have been little studied before. The peculiarity of the current period of
international relations is in that the sanctions are applied by different agents iteratively, in packages and stages. They have a dramatic effect on international relations. However, there are still no simple mathematical models describing the main effect of sanctions and allowing one to derive stability conditions.

Below we propose a simple model of repeated  sanctions and counter-sanctions mutual dynamics, motivated by and having similarities with the Osipov-Lanchester model. Stability conditions for the overall  systems are derived both for deterministic case and in presence of stochastic uncertainty.

\section{Deterministic model of the repeated sanctions dynamics}

Recall that the simplest Osipov-Lanchester \footnote{The model below was proposed independently during World War I by Mikhail Osipov in 1915 in Russia and by Frederick Lanchester in 1916 in Great Britain, see \cite{Helmbold93}} is a system of linear differential equations
$$
\dot R(t) = -\alpha G,
$$
$$
\dot G=-\beta R ,
\eqno(1)
$$
where  $R,G$ are numbers of the units representing power of each opponent, $\alpha, \beta$ are firing intensities (efficiencies). Using Euler’s method a discrete time version of (1) can be written as follows
$$
R_{n+1}=R_{n}-\alpha\Delta t G_n,
$$
$$
G_{n+1}=G_{n}-\beta\Delta t R_n.
\eqno(2)
$$
where $\Delta t$ is sampling interval. The system (2)  can be considered as a model for warfare in their own right, with each iteration corresponding to a separate battle.

The above equations are not suitable for description of economical or political sanctions since efficiency of a sanction package is mainly determined by its negative influence on the opponent with respect to the opponent's counter-sanction package at the previous stage. Given a diplomatic tradition
of ''mirror response'', the sanction intensity at the next stage should be almost similar or at least not weaker than the intensity of counter-sanction at the previous stage. Note that the meaning of the term ''not weaker`` depends on the public opinion stronger than on the real economic effect. This claim is illustrated with the whole international sanction story of the last decade and especially of the first half of 2022. In order to jump to the equations recall that strength (weakness) of the sanction depends mostly on the increase (decrease) of the sanction pressure during the previous sampling interval. Assuming linearity of all dependencies we arrive at the following equations:
$$
x_{n+1}=x_{n}+\alpha (y_n-y_{n-1})
$$
$$
y_{n+1}=y_{n}+\beta (x_n-x_{n-1}),
$$
$$
n=1,2,...
\eqno(3)
$$
where $x_{n},y_n$ are values of sanction and counter-sanction pressure at $n^{th}$ sampling instant,
$\alpha, \beta$ are positive  cross-gain parameters. Thus, equations (3) constitute a simple linear model of sanction dynamics. Cross-gains $\alpha, \beta$ in general differ from the corresponding coefficients in (1),(2). The question: how to choose $\alpha, \beta$ will be discussed later.

In order to transform (3) to a more convenient form
denote $v_n=x_{n}-x_{n-1},w_n=y_{n}-y_{n-1}$. Then (3) takes form of vector power series
$$
v_{n+1}=\alpha w_n,
$$
$$
w_{n+1}=\beta v_n,
$$
$$
n=1,2,...
\eqno(4)
$$
Shift $n$ to $n-1$ in the second equation of (4) and substitute it into the first one. Then perform similar procedure  with the first equation of (4). We come up with the following equations:
$$
 v_{n+1}=q v_{n-1},
 $$
 $$
 w_{n+1}=q w_{n-1},
 $$
 $$
n=2,3,...
\eqno(5)
$$
where $q=\alpha\beta>0$ is total gain of the system.

Analysis of the linear system  (5) yields  first yet important conclusions concerning asymptotic behavior of the process. The only parameter influencing stability of (5) is the gain $q$. The system is asymptotically stable, i.e. $v_n\to 0, w_n\to 0$ as $n\to\infty$ for any initial conditions $v_0,w_0$ if and only if $q<1$. If $q>1$ then the solutions tend to infinity for all nonzero initial conditions. The boundary case $q=1$ corresponding to Lyapunov stability has zero measure in the space of the parameters $\alpha,\beta$ and and can be neglected in practice.

The first important observation is that strategy parameters $\alpha,\beta$ of the opponents enter expression for $q$ symmetrically. It means that both opponents are equally responsible for stability of the process.
Second important observation is that each opponent is able to ensure stability of the process by means of proper choice of its strategy (decreasing its own gain $\alpha$ or $\beta$). In other words it is profitable for each agent to decrease $q$ since increase of $q$ leads to instability which is disadvantageous both for that agent and for the whole system.

Assume that $x_0,x_1$ and $y_0,y_1$ are given. Then solutions to (4) are represented explicitly as power functions
$$
 x_{2n}=(q^{n}-1)/(q-1)(y_1-y_0)/\beta,
 $$
 $$
 x_{2n+1}=(q^{n}-1)/(q-1)(x_1-x_0)/\beta
$$
$$
 n=1,2,...
 \eqno(6)
 $$
Relations (6) and similar expressions for $y_n$ represent long term cumulative effect of sanctions.
It seems however that sanctions of the past weakly influence  economic wealth and public opinion in long term since each opponent takes all measures to suppress the effect of the sanctions as soon as possible.

\section{Stochastic model of the repeated sanctions dynamics}

Idealized model (3) has an apparent drawback that it may be supersensitive to the parameter values $\alpha$ and $\beta$. In reality the values of $\alpha$ and $\beta$ depend in a complex way on many factors, including economical, political ones and on public opinion. Moreover the values
of $\alpha$ and $\beta$ may change from stage to stage. The simplest model for uncertainty is randomness. Therefore a model with random parameters is proposed and analyzed in this section.

Suppose that $\alpha$ and $\beta$ in (3) are replaced with
$\alpha+\xi_n$ and $\beta+\eta_n$, respectively. The stochastic version of (3) is as follows:
$$
v_{n+1}=(\alpha+\xi_n) w_n,
$$
$$
w_{n+1}=(\beta+\eta_n) v_n,
$$
$$
n=1,2,...
\eqno(7)
$$
The random errors are assumed to have zero means and bounded variances:
$$
E\xi_n=E\eta_n=0, E\xi_n^2=\sigma_x^2, E\eta_n^2=\sigma_y^2,
$$
$$
n=1,2,...
$$
where $E$ stands for mathematical expectation. Assume additionally that the errors of different opponents are uncorrelated: $E\xi_n\eta_m=0$ for all $m,n=0,1,2,...$ and analyze mean square stability of (6).

Under imposed assumptions the mean squares of the variables in (7) satisfy the relations
$Ev_{n+1}^2=E[(\alpha+\xi_n)(\beta+\eta_n)]^2Ev_{n-1}^2$. It means that the stability of the system depends on the averaged gain
$\bar q=E[(\alpha+\xi_n)(\beta+\eta_n)]^2$. Evaluation of the averaged gain yields:
$$\bar q=E(\alpha^2+2\alpha\xi_n+\xi_n^2)(\beta^2+2\beta\xi_n+\xi_n^2)$$
Taking into account zero correlations we obtain
$$
\bar q=(\alpha^2+\sigma_x^2)(\beta_y^2+\sigma_y^2).
\eqno(8)
$$
It is seen that the same conclusions hold for stochastic case if the gain $q$ is replaced with the averaged gain $\bar q$. It is seen also that stability conditions for stochastic case are more strict since $\bar q>q$.
Note that stability may be achieved in principle only under condition
$(\sigma_x^2+\sigma_y^2)<1$.

\section{Conclusions}

The paper contains only conceptual exposition of the proposed model and its properties.
Such important issues as evaluation of the model parameters and using the model in practice as well as taking into account nonlinear effects are beyond this discussion. 
However some qualitative conclusions can be drawn even from this simplistic model.

The first important conclusion is that {\bf both opponents are equally responsible for stability of the process}. This conclusion holds  both in deterministic and in stochastic case.
The author believes it is very important for international stability studies. If one of the opponents increases strength of its sanctions unlimitedly without coordination with the other opponent then instability occurs for sure (in stochastic case in mean square and with probability one).
It means that international stability may be achieved only by coordinated efforts of all 
countries (agents). Coordination should take into account that the  values of the parameters $\alpha,\beta$ are different since they depend on historical, cultural, geopolitical reasons.

The second conclusion is that in deterministic case each opponent is able to ensure stability of the process by means of proper choice of its own strategy (decreasing its own gain). However in stochastic case it is not always possible. Stability may be achieved only under condition $\sigma_x^2+\sigma_y^2<1$, i.e. uncertainty matters.

Nevertheless it is still profitable for each opponent to decrease $q$ since increase of $q$ leads to instability which is disadvantageous both for that opponent and for the whole system.



\end{document}